# BOSONIC AND FERMIONIC SINGULARITIES IN DIFFEOLOGY

PATRICK IGLESIAS-ZEMMOUR

ABSTRACT. We explore the differential geometry of the quadrant $C_2 = [0, \infty[^2$, equipped with the subset diffeology of $\mathbf{R}^2$. We show a striking dichotomy between differential forms and symmetric tensors. While differential forms on $C_2$ are simply restrictions of smooth forms on $\mathbf{R}^2$ (a "fermionic" behavior where singularities are hidden), symmetric tensors exhibit a "bosonic" behavior where singularities accumulate. We prove a decomposition theorem identifying exactly the singular parts: they are purely axial. Surprisingly, the mixed interaction term is forced to be regular by the symmetries of the corner. Finally, we introduce the notion of *singular capacity* to quantify the order of singularity a tensor can support.

## 1. INTRODUCTION

Diffeology offers a change of paradigm regarding singularities. Instead of defining smoothness via charts to a fixed local model, diffeology tests the space via smooth parametrizations (plots). A geometric object is smooth if its pullback by any smooth plot is smooth.

In this note, we investigate how this definition interacts with the simplest singularities: boundaries and corners. We focus on two spaces: the half-line $\Delta = [0, \infty[$ and the quadrant $C_2 = [0, \infty[^2$, both embedded in their respective Euclidean spaces.

We find that the algebraic nature of the geometric object determines its ability to "see" the singularity. We use the terminology "Fermionic" and "Bosonic" to highlight the contrast between the antisymmetric nature of forms (akin to the

*Date*: December 2025.
2020 *Mathematics Subject Classification*. Primary 58A40; Secondary 58A10, 58A32, 53A45, 58K70.
*Key words and phrases*. Diffeology, Differential spaces, Manifolds with corners, Differential forms, Symmetric tensors, Singularities, Singular capacity.
The author thanks the Hebrew University of Jerusalem, Israel, for its continuous academic support. He is also grateful for the stimulating discussions and assistance provided by the AI assistant Gemini (Google) during the preparation of this manuscript, and lastly to DeepSeek for fine tuning and final rendering.

Pauli exclusion principle) and the symmetric nature of tensors (akin to Bose-Einstein accumulation).

(1) **Fermionic behavior:** Differential forms are blind to the boundary. The antisymmetry forces them to vanish on singular strata in a way that makes them indistinguishable from restrictions of ambient forms.
(2) **Bosonic behavior:** Symmetric tensors detect the boundary. The symmetry allows singular terms to survive. However, we show that this accumulation is strictly limited to the diagonal terms of the tensor.

The paper is organized as follows. Section 2 treats the half-line $\Delta$, establishing the basic decomposition for symmetric 2-tensors and contrasting it with differential forms. Section 3 recalls the diffeological definition of a Riemannian metric and shows why singular tensors are excluded. Section 4 generalizes the decomposition to the quadrant $C_2$, revealing the regularity of the mixed term. Finally, Section 5 introduces the concept of *singular capacity*, providing a scaling argument that explains the observed hierarchy.

## 2. THE CASE OF THE HALF-LINE

Let $\Delta = [0, \infty[ \subset \mathbf{R}$ equipped with the subset diffeology, and let $x$ denote the standard coordinate.

**Proposition 2.1** (Differential forms on $\Delta$). *The inclusion $j \colon \Delta \hookrightarrow \mathbf{R}$ induces a surjection on differential forms. Any $1$-form $\alpha$ on $\Delta$ is of the form $\alpha = f(x)dx$ where $f$ is the restriction of a smooth function on $\mathbf{R}$. In particular, differential forms do not "see" the boundary.*

*Proof.* See [GIZ16]. The key observation is that if a plot $P(t)$ touches the boundary at $t_0$ (i.e., $P(t_0) = 0$), then $P'(t_0) = 0$ because P is non-negative. Hence the pullback $P^*(\alpha) = f(P(t))P'(t)dt$ vanishes to first order, hiding any potential singularity. □

The situation for symmetric tensors is different. Consider the symmetric 2-tensor

$$\tau_{\text{sing}} = \frac{dx \otimes dx}{x}. \tag{1}$$

**Lemma 2.2.** *The tensor $\tau_{\text{sing}}$ defined in (1) is a smooth diffeological symmetric $2$-tensor on $\Delta$.*

*Proof.* Let $P \colon U \to \Delta$ be a plot. We must verify that $P^*(\tau_{\text{sing}}) = (P'(t)^2/P(t))dt \otimes dt$ is smooth. If $P(t_0) > 0$, smoothness is trivial. Assume $P(t_0) = 0$. The following classical inequality (Glaeser–Landau) provides the required control. □

**Lemma 2.3** (Glaeser–Landau inequality). *Let $f : \mathbf{R} \to \mathbf{R}$ be a smooth non-negative function. For any compact interval $I$, there exists a constant $C = \sup_I |f''|$ such that for all $t \in I$:*

$$f'(t)^2 \leq 2C f(t).$$

*Proof.* Fix $t \in I$. By Taylor's theorem with Lagrange remainder, for any $h$,

$$f(t+h) = f(t) + h f'(t) + \frac{h^2}{2} f''(\xi),$$

for some $\xi$ between $t$ and $t + h$. Since $f(t+h) \geq 0$, the discriminant of this quadratic polynomial in $h$ must be non-positive:

$$f'(t)^2 - 2f(t) f''(\xi) \leq 0.$$

Taking the supremum of $|f''|$ on $I$ yields the claim. □

Applying this to $P$ (restricted to a compact neighborhood of $t_0$), we obtain

$$\left| \frac{P'(t)^2}{P(t)} \right| \leq 2 \sup |P''|,$$

which remains bounded. Moreover, if $P$ is flat at $t_0$ (i.e., $P'(t_0) = P''(t_0) = \cdots = 0$), then the ratio vanishes at $t_0$ together with all its derivatives. Hence $P^*(\tau_{\text{sing}})$ is smooth.

We now prove that $\tau_{\text{sing}}$ and regular tensors exhaust all possibilities.

**Theorem 2.4** (Decomposition on the half-line). *The space of differential symmetric $2$-tensors on $\Delta$ decomposes uniquely as*

$$\mathscr{S}^2(\Delta) = \mathbf{R} \cdot \tau_{\text{sing}} \oplus j^*\big(\mathscr{S}^2(\mathbf{R})\big).$$

*That is, every tensor $\tau \in \mathscr{S}^2(\Delta)$ can be written uniquely as $\tau = c \tau_{\text{sing}} + \tau_{\text{reg}}$, where $c \in \mathbf{R}$ and $\tau_{\text{reg}}$ is the restriction of a smooth symmetric $2$-tensor on $\mathbf{R}$.*

*Proof.* On the interior $]0, \infty[$, write $\tau = f(x) dx \otimes dx$ for some smooth $f : ]0, \infty[ \to \mathbf{R}$. Consider the square map $\text{sq} : \mathbf{R} \to \Delta$, $\text{sq}(t) = t^2$, which is a subduction. The pullback $\text{sq}^*(\tau)$ must be a smooth symmetric 2-tensor on $\mathbf{R}$:

$$\text{sq}^*(\tau) = f(t^2) d(t^2) \otimes d(t^2) = 4t^2 f(t^2) dt \otimes dt.$$

Define $g(t) = 4t^2 f(t^2)$. Then $g \in C^\infty(\mathbf{R})$ and, because sq is even, $g$ is even. By Whitney's theorem on even functions [Whi43], there exists $h \in C^\infty(\mathbf{R})$ such that $g(t) = h(t^2)$. Restricting to $t \geq 0$ and setting $x = t^2$, we have $h(x) = 4x f(x)$ for $x > 0$.

Expand $h$ in Taylor series at 0:

$$h(x) = h(0) + x k(x),$$

with $k \in C^\infty(\mathbf{R})$. Consequently,

$$4xf(x) = h(0) + xk(x) \implies f(x) = \frac{h(0)}{4x} + \frac{k(x)}{4}.$$

Thus on $]0, \infty[$,

$$\tau = \frac{h(0)}{4} \frac{dx \otimes dx}{x} + \frac{k(x)}{4} dx \otimes dx.$$

The first term is a multiple of $\tau_{\text{sing}}$; the second extends smoothly to $\mathbf{R}$. Uniqueness follows because $\tau_{\text{sing}}$ is unbounded at 0 whereas any element of $j^* \mathscr{S}^2(\mathbf{R})$ is bounded. $\square$

## 3. RIEMANNIAN METRICS

We recall the diffeological definition of a Riemannian metric [PIZ25, p. 224].

**Definition 3.1.** *A* Riemannian metric *on a diffeological space* X *is a smooth symmetric covariant 2-tensor g that is:*

(1) **Positive:** *For every $x \in X$ and every (germ of) path $\gamma \in \text{Paths}(X)$ with $\gamma(0) = x$, $g(\gamma)_0(1,1) \geq 0$.*
(2) **Definite:** *For every $x \in X$ and every plot $\gamma$ centered at $x$,*

$$g(\gamma)_0(1,1) = 0 \iff \forall \alpha \in \Lambda^1_x(X), \alpha(\gamma)_0(1) = 0,$$

*where $\Lambda^1_x(X)$ denotes the space of pointed 1-forms at $x$.*

**Proposition 3.2.** *On $\Delta$, a Riemannian metric cannot contain a singular part; i.e., if $g = c\tau_{\text{sing}} + \tau_{\text{reg}}$ is a Riemannian metric, then $c = 0$.*

*Proof.* Consider the plot $\gamma(t) = t^2$, centered at 0. For any pointed 1-form $\alpha = f(x)dx$, we have $\alpha(\gamma)_0(1) = f(0) \cdot 0 = 0$. Hence the definiteness condition requires $g(\gamma)_0(1,1) = 0$. But

$$\gamma^*(\tau_{\text{sing}}) = \frac{(2t\,dt)^2}{t^2} = 4dt^2,$$

so $\tau_{\text{sing}}(\gamma)_0(1,1) = 4 \neq 0$. Therefore $c$ must be zero. $\square$

**Remark 3.3.** *This argument highlights a key difference between arbitrary symmetric tensors and Riemannian metrics: the latter's non-degeneracy forbids singularities of type $1/x$ at the boundary. The "bosonic accumulation" is geometrically possible, but metrically excluded.*

## 4. The Case of the Quadrant

Let $C_2 = [0, \infty[^2 \subset \mathbf{R}^2$ with coordinates $(x, y)$. This space models a corner. We study the structure of $\mathscr{S}^2(C_2)$.

**Theorem 4.1** (Decomposition on the quadrant). *Every differential symmetric 2-tensor $\tau$ on $C_2$ decomposes uniquely as*

$$\tau = \tau_{\text{reg}} + \frac{A(y)}{x} dx^2 + \frac{B(x)}{y} dy^2,$$

*where $\tau_{\text{reg}}$ is the restriction of a smooth symmetric 2-tensor on $\mathbf{R}^2$, and $A, B \in C^\infty([0, \infty[)$.*

*Proof.* On the interior $]0, \infty[^2$, write

$$\tau = \alpha(x, y) dx^2 + \beta(x, y) dy^2 + 2\gamma(x, y) dx dy.$$

Consider the square map $\text{sq}: \mathbf{R}^2 \to C_2$, $\text{sq}(u, v) = (u^2, v^2)$, which is a subduction. Pulling back $\tau$ gives a smooth tensor on $\mathbf{R}^2$:

$$\text{sq}^*(\tau) = \alpha(u^2, v^2)(2u du)^2 + \beta(u^2, v^2)(2v dv)^2 + 2\gamma(u^2, v^2)(2u du)(2v dv)$$
$$= 4u^2 \alpha(u^2, v^2) du^2 + 4v^2 \beta(u^2, v^2) dv^2 + 8uv \gamma(u^2, v^2) du dv.$$

Define

$$\tilde{\alpha}(u, v) = 4u^2 \alpha(u^2, v^2), \quad \tilde{\beta}(u, v) = 4v^2 \beta(u^2, v^2), \quad \tilde{\gamma}(u, v) = 4uv \gamma(u^2, v^2).$$

These functions are smooth on $\mathbf{R}^2$. The map $\text{sq}$ is invariant under the group $\Gamma = \mathbf{Z}_2 \times \mathbf{Z}_2$ acting by sign changes on $(u, v)$. Consequently $\text{sq}^*(\tau)$ must be $\Gamma$-invariant.

**1. Axial terms.** The term $du^2$ is invariant under $(u, v) \mapsto (-u, v)$ and $(u, v) \mapsto (u, -v)$. Hence $\tilde{\alpha}$ must be even in both $u$ and $v$. By Whitney's theorem, there exists a smooth $K_\alpha$ on $\mathbf{R}^2$ such that $\tilde{\alpha}(u, v) = K_\alpha(u^2, v^2)$. Setting $x = u^2$, $y = v^2$, we have

$$4x \alpha(x, y) = K_\alpha(x, y).$$

Expand $K_\alpha$ in Taylor series with respect to $x$:

$$K_\alpha(x, y) = K_\alpha(0, y) + x Q_\alpha(x, y),$$

with $Q_\alpha$ smooth. Thus

$$\alpha(x, y) = \frac{K_\alpha(0, y)}{4x} + \frac{Q_\alpha(x, y)}{4}.$$

The first summand is of the form $A(y)/x$ with $A(y) = K_\alpha(0, y)/4$ smooth on $[0, \infty[$. The second summand is regular. The same reasoning applies to $\beta$, yielding a term $B(x)/y$.

**2. Cross term.** The product $du dv$ changes sign under $(u, v) \mapsto (-u, v)$ and under $(u, v) \mapsto (u, -v)$. Hence $\tilde{\gamma}$ must be odd in both $u$ and $v$. A smooth function with this symmetry can be written as $\tilde{\gamma}(u,v) = uv\,K_\gamma(u^2, v^2)$ for some smooth $K_\gamma$. Comparing with the definition of $\tilde{\gamma}$,

$$4uv\gamma(u^2, v^2) = uv\,K_\gamma(u^2, v^2),$$

so for $uv \neq 0$ we have $4\gamma(x, y) = K_\gamma(x, y)$. Since $K_\gamma$ is smooth, $\gamma$ extends to a smooth function on all of $C_2$. Therefore the cross term $\gamma(x,y)dx\,dy$ is necessarily regular.

The uniqueness follows from the observation that the singular parts $A(y)/x$ and $B(x)/y$ are uniquely determined by the behavior of $\tau$ near the axes. $\square$

**Remark 4.2.** *The theorem reveals a remarkable* selection rule*: singularities can only appear in the diagonal components $dx^2$ and $dy^2$, while the off-diagonal term $dx\,dy$ is forced to be regular. This is a direct consequence of the different transformation properties under the sign-change symmetries of the corner.*

## 5. SINGULAR CAPACITY

The previous results suggest a hierarchy: the degree of a tensor determines the strength of the singularity it can support. This can be understood through a scaling argument.

Let $P(t)$ be a plot of $\Delta$ with $P(t_0) = 0$. Because $P(t) \geq 0$, the derivative $P'(t)$ must vanish at $t_0$; in fact, $P(t)$ behaves like $(t - t_0)^{2m}$ for some $m \geq 1$ (typically $m = 1$). Hence $P'(t) \sim (t - t_0)^{2m-1}$, and the pullback of a $k$-tensor involves a factor $(P'(t))^k \sim (t - t_0)^{k(2m-1)}$. On the other hand, a singularity of type $1/x^p$ corresponds to $1/P(t)^p \sim 1/(t - t_0)^{2mp}$. For the pullback to remain smooth, we need the vanishing order of the numerator to outweigh the divergence of the denominator:

$$k(2m-1) \geq 2mp \quad \text{for all } m \geq 1.$$

The worst case is $m = 1$, giving $k \geq 2p$. Thus the maximal admissible pole order is $p = \lfloor k/2 \rfloor$.

**Definition 5.1** (Singular capacity). *The* singular capacity *of a covariant $k$-tensor on the half-line $\Delta$ (or more generally on a boundary of codimension 1) is the integer*

$$\varkappa(k) = \left\lfloor \frac{k}{2} \right\rfloor,$$

*representing the maximum integer $p$ such that a singularity of type $1/x^p$ can be supported while maintaining smoothness in the diffeological sense.*

**Example 5.2.**

- $k = 0$ *(functions) and $k = 1$ (1-forms): $\varkappa = 0$; no poles allowed.*
- $k = 2$ *(symmetric 2-tensors): $\varkappa = 1$; a simple pole $1/x$ is allowed (as in $\tau_{\text{sing}}$).*

- *$k = 3$: $\varkappa = 1$; again a simple pole is allowed, but not $1/x^2$.*
- *$k = 4$: $\varkappa = 2$; a pole of order 2 becomes possible.*

**Remark 5.3.** *This notion provides a diffeological analogue of distribution theory: the plots act as test curves, and the vanishing of derivatives at the boundary imposes natural regularization. The singular capacity quantifies how much "roughness" a tensor of a given rank can accommodate.*

## CONCLUSION

The quadrant $C_2$ exhibits a rich diffeological structure. While differential forms behave "fermionically" (blind to the corner), symmetric tensors behave "bosonically" but with a selection rule. The diagonal components of a symmetric tensor can blow up at the boundary (creating a singular geometry that opens up like a horn), but the off-diagonal terms must remain regular. This suggests that Diffeology imposes a specific orthogonality constraint on singular tensors at the corner.

The concept of singular capacity offers a unifying perspective, linking the algebraic rank of a tensor to the order of singularities it can support. This hierarchy may have implications for the study of singular Riemannian metrics, stratified spaces, and quantization in the presence of boundaries.

EINSTEIN INSTITUTE OF MATHEMATICS,, THE HEBREW UNIVERSITY OF JERUSALEM,, CAMPUS GIVAT RAM,, 9190401 ISRAEL.

*Email address*: piz@math.huji.ac.il

*URL*: http://math.huji.ac.il/~piz